\documentclass[11pt]{amsart}
\usepackage{latexsym,amssymb,amsmath}
\def\oio{\o_{i_0}}\def\aio{\a_{i_0}}

\def\t{\tau}
\def\s{\sigma}

\def\ep{\epsilon}

\def\ee#1{e_{#1}}

\def\trank{\text{rank}}

\def\BC{\mathbb C}\def\BF{\mathbb F}\def\BO{\mathbb O}\def\BS{\mathbb S}

\def\BP{\mathbb P}
\def\pp#1{\mathbb P^{#1}}
\def\fa{\mathfrak a}
\def\fb{\mathfrak b}

\def\pp#1{{\mathbb P}^{#1}}
\def\tdim{\rm dim}
\def\hd{,...,}
\def\ww{\wedge}

\def\inv{{}^{-1}}

\def\cB{{\mathcal B}}\def\cA{{\mathcal A}}

\def\cQ{{\mathcal Q}}

\def\cR{{\mathcal R}}\def\cM{{\mathcal M}}
\def\cS{{\mathcal S}}

\def\cO{{\mathcal O}}
\def\CC{\mathbb C}

\def\BZ{\mathbb Z}
\def\SS{\mathbb S}

\def\11{\mathbf 1}

\def\QQ{\mathbb Q}

\def\fsl{{\mathfrak {sl}}}

\def\ff{{\mathfrak f}}

\def\fg{{\mathfrak g}}
\def\fn{{\mathfrak n}}
\def\fp{{\mathfrak p}}

\def\ft{{\mathfrak t}}
\def\fl{{\mathfrak l}}
\def\l{\lambda}
\def\a{\alpha}

\def\o{\omega}

\def\b{\beta} 
\def\g{\gamma}
\def\s{\sigma}

\def\n{\nu}
\def\d{\delta}

\def\ot{{\mathord{\,\otimes }\,}}
\def\op{{\mathord{\,\oplus }\,}}
\def\otc{{\mathord{\otimes\cdots\otimes}\;}}

\def\ra{{\mathord{\;\rightarrow\;}}}

\def\La#1{\Lambda^{#1}}
\def\tim{\text{Image}\,}

\def\tdim{\text{dim}\,}

\def\tker{\text{ker}\,}

\def\tmod{\text{ mod }}
\def\tmin{\text{ min }}
\def\tmax{\text{ max }}

\def\trank{\text{rank}\,}

\newtheorem{theorem}{Theorem}[section]
\newtheorem{proposition}[theorem]{Proposition}
\newtheorem{lemma}[theorem]{Lemma}
\newtheorem{corollary}[theorem]{Corollary}
\newtheorem{conjecture}[theorem]{Conjecture}

\theoremstyle{definition}

\newtheorem{example}[theorem]{Example}

\theoremstyle{remark}
\newtheorem{remark}[theorem]{Remark}

\def\ctimes{\times\cdots\times}

\begin{document}

\title{On tangential varieties of
rational homogeneous varieties}
\author{J.M. Landsberg and Jerzy Weyman}
\begin{abstract} We determine which tangential varieties
of homogeneously embedded rational homogeneous varieties
are spherical. We determine the
homogeneous coordinate rings
and rings of covariants of the tangential
varieties of   homogenously embedded compact Hermitian
symmetric spaces (CHSS). We give bounds on the degrees of
generators of the ideals of tangential varieties of
CHSS and obtain more explicit infomation about
the ideals in certain cases.
\end{abstract}

 \thanks{Supported respectively by NSF grants DMS-0305829 and DMS-0600229}
\email{jml@math.tamu.edu, j.weyman@neu.edu}

\maketitle

\section{Introduction}

Let $K$ be an algebraically closed field of characteristic zero, let $G$ be a semisimple algebraic group over $K$, and let
$P$ be a parabolic subgroup. Consider the homogeneous space $X=G/P$ embedded equivariantly
as the orbit of a highest weight   line
 in a projective space $\BP  V$, where $V$ is
an irreducible $G$-module.

We investigate   properties of 
the {\it tangential variety} $\t (X)\subset \BP V$ which is, 
by definition, the union of the points
on  the embedded tangent lines (i.e. $\BP^1$'s) to $X$.

In section  \ref{sphericsect},   we determine when $\t (X)$ is spherical
(Theorem \ref{spherethm}).   We show $\t (X)$ is 
$G$-spherical iff $X$ admits the structure
of a compact Hermitian symmetric space (CHSS)  (possibly
for a larger group $G'\supseteq G$) except for the case 
that $X$ contains $ G_2/P_1 $ as a factor.
The results of this section are proved using   properties of the
projective second fundamental form and combinatorics of root
systems, with the exception of the case where $G$ is of type
$A_n$, where we use an additional argument.

For the remainder of the paper we restrict to
the case where $V$ is a {\it generalized cominuscule}
$G$-module. One can define
a generalized cominuscule module  $V$
by the property that for $G/P=X\subset \BP V$,
the tangent bundle
$TX$ is an irreducible homogeneous vector bundle. 
With this definition, a {\it cominuscule module} is
a generalized cominuscule module where
moreover $G$ is
simple and $V$ is a fundamental representation.
We refer to $X\subset\BP V$ as a {\it (generalized) cominuscule
variety}.
A generalized cominuscule variety admits the
structure of a $G$-compact Hermitian symmetric
space, and when we refer to the {\it rank}
of a generalized cominuscule variety, we mean
the rank of the corresponding CHSS. Note that
the rank can be defined purely in terms of the projective
structre as the length of the osculating sequence of $X$,
see \cite{LM0}.

In  sections \ref{irredchsssect},\ref{multsect}
and \ref{reduciblecominsect},
assuming $X$ is generalized cominuscule,   we study the coordinate rings of its tangential variety $\t (X)$. We use the fact that in the case when $\t (X)$ is   non-degenerate, the cone over it has a natural desingularization which is the total space of homogeneous vector bundle. We use    methods from
\cite{weyman}, which reduce  the calculation of syzygies 
 to 
the calculation of   sheaf cohomology
groups of certain vector bundles. We calculate some of this cohomology to determine the   the coordinate ring of   $\t (X)$  as a $G$-module (Theorems \ref{irrcoordgens}, \ref{redcoordgens}), which
turns out to be uniform over cominuscule varieties
of the same rank.  

In  section \ref{degensect}  we apply results from \S\ref{irredchsssect},\S\ref{multsect}
and \S\ref{reduciblecominsect}, combined with    deformation results of Grosshans 
\cite{Gr} to deduce  upper bounds on the degrees of generators of the defining ideals of $\tau (X)$ when $X$ is generalized cominuscule.

  Finally in section \ref{examsect}  we draw   consequences from the previous
sections to give  results on the defining ideals,
describing certain cases explicitly and bounding the
degrees of the generators of the ideals in all cases.

We describe the content in more detail below.

\subsection{Notation}\label{notation}
For a variety $Z\subset \BP V$, we let
$\hat Z\subset V$ denote the corresponding cone.
For $z\in Z$ a smooth point, we let
$\hat T_zZ\subset V$ denote the affine tangent space,
$T_zZ$ the Zariski tangent space, and
$\tilde T_zZ=\BP(\hat T_zZ)\subset \BP V$ the embedded
tangent projective space. The two are related
by $T_zZ=(\hat T_zZ/\hat z)\ot \hat z^*$.
Let $K$, $G$, $P$, $X$ be as in the first paragraph.
We use German letters to denote Lie algebras
associated to algebraic groups.
We use the ordering of roots as in \cite{bour}.
   The fundamental weights
and  the simple roots  of $\fg$ are respectively
denoted  $\o_i$  and
$\a_i$.
  $P_k$ denotes the maximal parabolic of $G$ obtained by
deleting 
the root spaces corresponding to  negative   roots having a nonzero coefficient on the simple
root $\a_k$.
More generally,
for $J=(j_1\hd j_s)$,   $P_J$ denotes the parabolic obtained
by deleting the negative root spaces     having a nonzero coefficient on any of the simple
roots
$\alpha_{j_1}\hd \alpha_{j_s}$.
$\Lambda_{\fg},\Lambda_{G}$
respectively denote the weight lattices of $\fg$, $G$,
and $\Lambda^+_{\fg}\subset \Lambda_{\fg}$, $ \Lambda^+_G\subset \Lambda_G$ the dominant weights. We let $L\subset P$ be a (reductive)
Levi factor and $\ff=[\fl,\fl]$ a semi-simple Levi
factor.  We write $\fp=\fl +\fn$, where $\fn$
is nilpotent.  

When dealing with $\fa_n$-modules we sometimes use partions
to index highest weights, with the dictionary
$\pi=(p_1\hd p_{n+1})$ corresponds to the
weight $(p_1-p_2)\o_1+(p_2-p_3)\o_2+\cdots + (p_n-p_{n+1})\o_n$.

\subsection{Sphericality}  Recall that a
normal projective $G$-variety
$Z$ is {\it $G$-spherical} if some (and hence any) Borel subgroup
  $B$ of $G$ has a dense orbit in $Z$.
(We emphasize the group $G$
in our terminology because certain varieties we study
will be $G$-varieties for several different groups.)
 Equivalently, $Z$ is
spherical if for all degrees $d$, $K[Z]_d$,
the component of the coordinate ring of $Z$ in degree $d$, is a multiplicity
free $G$-module, see \cite{briont}.
Note that this property for $\t(X)$ {\it a priori} depends both
on $G$ and the embedding of $X$.

\begin{theorem}\label{spherethm} 
Let $X=G/P \subset \BP V$ be a 
homogeneously embedded rational homogeneous variety.
Then $\t (X)$ is $G$-spherical iff $X$ admits
the structure of a CHSS,
and no factor of $X$ is $G_2/P_1$. 
\end{theorem}

Note that the varieties $C_n/P_1$, $B_n/P_n$, which are not
cominuscule, have spherical tangential varieties.

\bigskip

\subsection{Ideals and singularities}\label{isingsect}

For any smooth variety $X\subset \BP V$, if the tangential
variety is {\it strongly nondegenerate} in the sense that a general
point of $\t (X)$ lies on a unique tangent line, then
$\t (X)$ admits a natural desingularization by the  projective bundle
$\tilde T  X $, whose fiber over $x\in X$ is the embedded tangent
projective space $\tilde T_xX\subset \BP V$. The associated
vector bundle $\hat T X $ is a subbundle of
the trivial bundle $V\ot \cO_X$. It is a desingularization of the affine cone $ \hat\t( X)\subset V$.
In the case
$X$ is homogenous and homogeneously embedded, this desingularization
is  an example of  
what Kempf \cite{kempf} called \lq\lq the collapsing of a
homogeneous vector bundle\rq\rq . We apply
the methods of \cite{weyman} to study the cone ${{\hat\t}(X)}$ via this
desingularization in the cases of generalized
cominuscule varieties. The essential point is
that the minimal free resolution
of the ideal of $\hat T X  \subset V\ot \cO_X$ is given by a
Koszul complex, and we can \lq\lq push down\rq\rq \ this
information to get information about the ideal of $\t(X)$.

\begin{proposition}
Let $X=G/P\subset \BP V$ be a homogeneously
embedded rational homogeneous variety. Then $\t(X)$ is strongly nondegenerate
except in the case when $X$ is generalized cominuscule, has rank  two
and the embedding is the minimal homogeneous embedding,
or $X=\pp n$ and the embedding is minimal or quadratic
Veronese.
\end{proposition}

This proposition is essentially \lq\lq known to the
experts\rq\rq\ because the first candidates for
non-strongly degenerate $\t(X)$, namely the
rank $3$ cominuscule varieties appearing in the third row
of the Freudenthal magic chart, are known
to be strongly nondegenerate (in fact they are \lq\lq
one apparent double point\rq\rq). In any case it
is easy to verify on a case by case basis.

In other words, for 
 cominuscule varieties $X\subset \BP V$, $\t(X)$ is strongly nondegenerate
except when it coincides with the secant variety of
$X$, $\s(X)$, which is the Zariski closure of
all points on all secant lines to $X$.  The
cominuscule rank two case is well
understood, see \S\ref{scorza}.

Recall that a variety $Y$ 
over a field of characteristic zero has {\it rational singularities}
if it is normal
and it admits a desingularization $\pi : Z\ra Y$ such that
$R^i\pi_*\cO_Z=0$ for $i>0$.

\begin{theorem}\label{tauthm}

Let $X\subset \BP V$ be a  
rank $r\geq 3$ cominuscule variety. Then
\begin{enumerate}

\item  ${ \t}(X) $ is normal, with rational singularities.

\item
The coordinate ring $K[\t (X)]$ has a uniform decomposition into
irreducible modules given by theorem \eqref{irrcoordgens}.1 in the irreducible
case and  theorem  \eqref{redcoordgens}.1 in the reducible case.

\item The ring of covariants of the coordinate ring $K[\t(X)]$  has generators described by 
 theorem \eqref{irrcoordgens}.2 in the irreducible
case and  theorem  \eqref{redcoordgens}.2 in the reducible case.

\item If $X$ is irreducible, the ideal of $\t (X)$ is generated
 in degrees at most $4r-4$ and if it is reducible with its highest rank
factor having rank $r_0$, then the ideal of $\t(X)$ is generated in
degrees at most $\tmax (6,4r_0)$.

\end{enumerate}

\end{theorem}

Our uniform degree bound  follows from adapting a deformation argument applicable to    spherical varieties. This idea dates back to 
unpublished work of Luna and to Brion's thesis
\cite{briont}. We use the version of Grosshans
\cite{Gr}. This is explained in \S\ref{degensect}.

In \S\ref{multsect} we describe the generators of rings of covariants for the tangential varieties of multiple embeddings of $X$.

We     give more precise information
about the ideals and singularities in  some cases in \S\ref{exsect}.

\medskip

We summarize from \cite{weyman} (5.1.1-3,5.4.1) the
results we will use:
\begin{theorem}\cite{weyman}\label{weymanthm}
 Let $Y\subset \BP V$ be a
variety and suppose there is a
projective variety $B$ and a vector bundle
$E\ra B$ that is a subbundle of a trivial bundle
$\underline V \ra B$ with $\underline V_z\simeq V$ for
$z\in B$ such that $  \BP E\ra   Y$ is a desingularization
of $Y$.  Write $\eta=E^*$ and $\xi=(\underline V/E)^*$.

If
the sheaf cohomology groups
$H^i(B,S^d\eta)$ are all zero for $i>0$
and if the linear maps
$H^0(B,S^d\eta)\ot V^*\ra H^0(B,S^{d+1}\eta)$
are surjective for all $d\geq 0$, then
\begin{enumerate}
\item
$\hat Y$ is normal, with rational singularities 

\item The coordinate ring $K[\hat Y]$ satisfies
to $K[\hat Y]_d\simeq H^0(B,S^d\eta)$.

\item The 
vector space of minimal generators of the ideal of $Y$ in
degree $d$ is isomorphic to  
$H^{d}(B,\La{d+1}\xi)$.

\item If moreover $Y$ is a $G$-variety
and the desingularization is $G$-equivariant,
then the identifications above are as $G$-modules.
\end{enumerate}
\end{theorem}

In our situation
$\xi:= (V\ot \cO_X/\hat TX)^* $   and
$\eta:=(\hat TX)^*$. The  bundles $\xi$ and $\eta$ are homogeneous
and indecomposable   but not irreducible,
so we  first calculate 
the cohomology of the corresponding graded bundles
using the Bott-Borel-Weil theorem and then pass
to the cohomology we are interested in using
methods from \cite{Orub}.

\medskip

\subsubsection*{Acknowledgments}
We   thank Michel Brion for pointing out the crucial reference \cite{Gr}, and   Laurent Manivel, Dmitry Panyushev
and the anonymous  referee for  very useful remarks.

\section{Determination when $\t (X)$ is spherical}\label{sphericsect}

In this section we reduce the calculation of whether or not
$\t (X)$ is spherical to a calculation if a linear representation
is spherical, i.e., {\it visible} in the language of
\cite{Kac}. Then we use the list in \cite{Kac} to determine the
cases where $P$ is maximal, and we determine   other cases
when $P$ is not maximal by hand.
 
For any smooth variety $X\subset \BP V$, let $x\in X$ be a general point and let $w\in T_xX$ be a generic
vector.
Then $\tdim \t(X)= \tdim X+\trank II_w$, where $II\in S^2T^*_xX\ot N_xX$ is the projective
second fundamental form of $X$ at $x$ and $II_w: T_xX\ra N_xX$ is the map $u\mapsto II(u,w)$,
see \cite{GH} or \cite{IvL}, Proposition 3.13.3.

\medskip

In our situation,  let $x=[v]$ correspond to a highest weight vector.  Then
$T_xX\simeq \fg.v \tmod v$ (up to a twist by a line bundle which we ignore  
throughout this section). Let $U_1, U_2\in \fg$, then   $II (U_1.v, U_2.v)= U_1.(U_2.v)\tmod \fg.v$,
see \cite{LM0}, Proposition 2.3.

Say the parabolic $\fp$ is obtained by deleting the root spaces of negative  roots having nonzero
coefficient on the simple roots $\b_1\hd \b_p$. These simple roots induce a $\BZ^p$ grading on
$\fg$ which induces multi-filtrations on $T_xX$ and $N_xX 
\simeq V/\hat T_xX$.

  Let $L\subset P$ be a Levi
factor, i.e., a maximal reductive subgroup of $P$, and
write $T=T_xX=\op_IT_I$ for the decomposition of $T$ as an $L$-module.
For $I=(i_1\hd i_p)$,   $T_I$ is the
sum of root spaces for roots $\g$ such that when $\g$ is expressed as sum of
simple roots that $i_j$ is the coefficient of $\b_j$ in the expression. 

\medskip

A $G$-variety $Z\subset \BP V$ is $G$-spherical iff there exists a Borel $B\subset G$ such that
there is an open $B$ orbit in $Z$, or equivalently, letting $z\in Z$ be a general point,
and $\overline z\in V$   a corresponding vector in the line $z\in \BP V$,
$Z$ is $G$-spherical 
iff $\fb.\overline z=\hat T_{ z}Z\subset V$, where $\hat T_zZ$ is the affine tangent space.  

Note that $\hat \t(G_2/P_1)=V_{\o_1}$ and $V_{\o_1}$ is
not  a visible (spherical) $G_2$-module by \cite{Kac}.

We choose the Borel $B$  containing
the negative roots. We compare the spaces $\hat T_{v+U.v}\t (X)$ and
$\fb.(v+U.v)$ inside $V$.
An arbitrary element of $\hat T_{v+U.v}\t (X)$ is of the form
$cv+U_1.v+U_2.U.v$ for $c\in \BC$,  $U_1 ,U_2\in \fg$ and without loss of generality we
may take $U_1 ,U_2\in \fg_{-}$.

Let $\fb_0\subset \fg_0$ denote the component of $\fb$ in $\fg_0$ so
$\fb=\fb_0+\sum_{I} \fg_{-I}$. For $b\in\fb$, write
$b=b_0+b_I$ with $b_0\in \fb_0$.

Similarly, write $U=\sum U_I$ with $U_I\in \fg_{-I}$, so
$$
U.v= \sum U_I.v
$$
with $U_I.v\in T_I$.
Now consider $U_2.U.v$.

Assume  
 that $\t(X)$ is of its expected dimension so $II_{U.v}$ is injective.
(The cases where $\t (X)$ is not of the expected dimension are
always spherical (excepting $G_2/P_1$), which will follow by re-embedding the variety in
such a way that $\t(X)$ is nondegenerate by the argument below.)
Then each $(U_2)_I$ must map injectively on $U.v$ and the vectors
$(U_2)_IU_J.v$   must be independent as elements of $N_xX$.

Thus each $b_J$   is used to fill $II_{U.v}(T)\subset N_xX$ exactly,
and the only question that remains is if the vectors
$$
b_0.U.v
$$
fill $\hat T_xX \tmod \hat x$. Since $b_0.v\subset \hat v$, we are
reduced to considering $[b_0,U]$.
We have proved:

\begin{lemma}\label{redlem}
Let $X=G/P$ be homogeneous.
 Assuming $\t (X)$ is nondegenerate, it is spherical iff for some  
$U\in \fg_{-}$, that $[\fb_0,U]=\fg_{-}$.
\end{lemma}

Lemma \ref{redlem} reduces the problem to a linear problem:  determining which $T_xX$,
considered as  $L$-modules, are $L$-spherical, where 
  $L\subset P$ is a Levi factor. The spherical irreducible $L$-modules were already determined
by Kac \cite{Kac}, theorem 3. Examining his list we immediately conclude
all cominuscule $X$ have $\t (X)$ spherical when $\t (X)$ is nondegenerate. But
 for the same $X=G/P$ in a smaller embedding,
this implies that $\t (X)$ is also spherical in the smaller embedding. Since the only
examples of homogeneous rational varieties with degenerate secant varieties occur
among the cominuscule varieties
 and   $G_2/P_1$, for all remaining cases we only need  determine if
$T_xX$   is $L$-spherical.
Recall the notation $P_I$ from \S\ref{notation}.

\begin{proposition} If $\t (G/P_{I})$ is not spherical, then $\t (G/P_J)$
is not spherical for any set $J$ of simple roots containing the
simple roots corresponding to $I$.
\end{proposition}

\begin{proof}
This is clear because $\fb_0^I\supset \fb_0^J$ but $\fg_{-}^I\subset \fg_{-}^J$.
\end{proof}

  Kac's list already eliminates all non-cominuscule $X$ except for:
$A_n/P$ (any $P$), $B_n/P_n$, $B_n/P_{1,n}$, $C_n/P$ (some $P$,
see below), $D_n/P_{1,n}$, $F_4/P_4$, $G_2/P_1$.

Among these  $F_4/P_4$ is immediately eliminated by dimension considerations as
 $\fl$ is  $\fb_3\op\BC$   and the $L$-module 
is $V_{\o_1}\op V_{\o_3}$, see \cite{LM0} .

\bigskip

\begin{proposition} $\t (B_n/P_n)$ is spherical.\end{proposition}
\begin{proof}
$T_x(B_n/P_n)\simeq \BC^n\op \La 2\BC^n$ and we must examine the action of
$\fb_0\supset \fb (\fsl_n)$ on it. Consider the vector $v=e_1\op e_1\ww e_2$, the
sum of two highest weight vectors. It is clear $\fb (\fsl_n).v$ is $T_x(B_n/P_n)$.
\end{proof}
 
\begin{proposition} The tangential varieties of $D_n/P_{1,n}$ and $B_n/P_{1,n}$
are not spherical.
\end{proposition}
\begin{proof}
For the $D_n/P_{1,n}$ case $T_xX$ splits into three $L$-modules. Letting
$V=V_{\o_1}=\BC^{n-2}$, they are 
$V \op V_{\o_{ n-2}}\op V_{\o_1+\o_{n-2}}$.
Let $(v,U,Z)\in T_xX$. In order to have $\fb_0$ cover $T_1$, we must
have $v=\ee 1$. Similarly to cover $T_{1,1}$ we must have 
$Z$ contain a summand of the form $e_1\ot (e_1\ww F)$. But now we see it
is impossible to have a vector of the form $(e_2,U, e_1\ww F)$ in 
$\fb_0.(v,U,Z)$. The $B_n$ case is similar.
\end{proof}

\bigskip

\begin{proposition} The only homogeneous $C_n$-varieties
having spherical tangential varieties are  $C_n/P_1, C_n/P_n$.
\end{proposition}

\begin{proof}
$\t(C_n/P_n)$, $\t(C_n/P_1)$  are
spherical by Kac's list (using $\fb_0$ in the first case
and $\fb$ in the second).  We will  rule out all
other  maximal parabolics, and
once having done so, the only other possibility would be $\t(C_n/P_{1,n})$, but
this is easily eliminated by dimension considerations.
(Note that Kac's list immediately implies $\t(C_n/P_k)$ is not
spherical if $k,n-k>4$.)

In terms of matrices, 
\begin{equation}
\fg_0= \left\{ \begin{pmatrix} a & 0 & 0 &0\\
0&b_1&0 &b_2\\ 
0&0&-{}^ta &0\\0&b_3&0& -{}^tb_1\end{pmatrix} \mid a\in \fsl_k,\ 
  b_1 \in\fsl_{n-k}, b_2={}^tb_2,\ b_3={}^tb_3
\right\}
\end{equation}
  
\begin{equation}
\fb_0= \left\{ \begin{pmatrix} a & 0 & 0 &0\\
0&b_1&0 &0\\ 
0&0&-{}^ta &0\\0&b_3&0& -{}^tb_1\end{pmatrix} \mid a\in \fb( \fsl_k),\ 
   b_1 \in \fb(\fsl_{n-k}),\ b_3={}^tb_3   
\right\}
\end{equation}
and   
\begin{equation}\label{Tsymbol}
T_{[v_{\o_k}]}(C_n/P_k)=   \begin{pmatrix} 0 & 0 & 0 &0\\
t&0&0 &0\\ 
T_2& {}^tt'&0 & {}^tt\\t'&0&0& 0\end{pmatrix}  
\end{equation}
with $t,t'\in Mat_{k\times n-k}$, $T_2\in S^2\BC^k$. 
Note that $T_1=t+t'$. 

Fixing an initial vector with components $(\t,\t',\t_2)$, the action of $\fb_0$ 
provides
\begin{align*}
 t&= b_1\t-\t a\\
T_2&= -( \t_2a+{}^t(\t_2a))\\
 t'&= -{}^tb_1\t'+b_3\t-\t'a
\end{align*}

Examining the $T_2$ term we see both $T_2$ and $a$ have dimension $\binom {k+1}2$ and
that we may use $a$ to exactly fill $T_2$. But then $t$ may be filled only by using
the lower diagonal matrix $b_1$, but   this is not possible 
when $k>1$.

\smallskip
 
Now we turn to the case of $A_n$.
Consider
  $A_{n-1}/P_{k,m}$. Let $k'=m-k$ and $k''=n-m$. Then $T$ may be thought of as the union of
three vector spaces consisting of the lower $k'\times k$, $k''\times k$
and $k''\times k'$ blocks. Note the naive bound for dimension reasons that
$$
kk'+kk''+k'k''\leq \binom{k+1}2 +  \binom{k'+1}2 + \binom{k''+1}2 -1.
$$
Use index ranges $1\leq i,j\leq k$, $k+1\leq s,t\leq m$
and $m+1\leq u,v\leq n$. Write the matrices to be filled as having elements
\begin{align}
 t^s_i &= \t^s_ja^j_i - a^s_t\t^t_i\label{one}\\
t^u_i &= \t^u_ja^j_i - a^u_v\t^v_i\label{two}\\
t^u_s &= \t^u_ta^t_s - a^u_v\t^v_s\label{three}
\end{align}
where the $\t$'s are given (generic) constants, and
given a set of $t$'s we want to determine if there
exists a set of $a$'s that produces them.
Here $a^A_B=0$ if $A<B$, $\sum_{A=1}^na^A_A=0$ and otherwise the entries 
$a^i_j,a^s_t,a^u_v$ of $\fb_0$ are independent. 
Now despite there being more unknowns than  equations to solve
in many cases, these equations are never compatible.
We illustrate with the adjoint case of $\fsl_{n+2}$,
we shift indices, having them run from $0$ to $n+1$,
so $1\leq s,t\leq n$.
We have
\begin{align}
 t^s_0 &= \t^s_0a^0_0 - a^s_t\t^t_0\label{one1}\\
t^{n+1}_{n+1} &= \t^{n+1}_0a^0_0 - a^{n+1}_{n+1}\t^{n+1}_0\label{two2}\\
t^{n+1}_s &= \t^{n+1}_ta^t_s - a^{n+1}_{n+1}\t^{n+1}_s\label{three3}
\end{align}
Label the equations
\eqref{one1},\eqref{two2},\eqref{three3} respectively by $(s,0)$, $(n+1,n+1)$
and $(0,s)$.
If we consider the following linear combination
of the right hand sides
$$
\sum_s t^{n+1}_s (s,0) -(\sum_st^{n+1}_st^s_{n+1})(n+1,n+1)
-\sum_st_{n+1}^s (0,s)
$$
we get zero, which shows there are choices of $t$'s for
which there does not exist a solution.
One can write out a proof of the general case
similarly, 
but we  instead include
a different proof, which, while using more machinery,
points out an explicit failure of sphericality.

\medskip

 (The following proof is best read after reading \S\ref{irredchsssect}.)

Consider $X=A_n/P_I\subset \BP V$
  for some $I=\lbrace i_1 ,\ldots ,i_s\rbrace\subset[1,n]$. 
Consider
  the desingularization of the affine cone $ \hat\t( X)\subset V$ given by the total space of a vector bundle $\hat T X $ constructed in  
\S\ref {isingsect}. Continuing the notation of \S\ref{isingsect}
with $\eta=(\hat TX)^*$,
$$
gr (\eta) = \cO_X(1)\oplus (\cO_X(1)\ot T^*X).
$$
By \cite{weyman} (5.1.2b) and (5.1.3.a),  the normalization of the ring of coordinate functions on $\t(X)$ is isomorphic to 
$H^0 (G/P ,Sym(\eta ))$.  Thus it is enough to show that for a non-maximal parabolic, that the algebra $H^0 (G/P ,Sym(\eta ))$ 
is not multiplicity free.

Recall   the notations $\o_i, \ep_i,\a_i$
from \cite{bour}:

\begin{proposition}
Assume that $s\ge 2$. The representation $V^*_{2(s\omega_{i_1}+(s-1)\omega_{i_2}+\ldots +\omega_{i_s})-\epsilon_{i_1}+\epsilon_{i_s+1}}$
occurs in $H^0 (G/P_I ,S^2 \eta )$ with multiplicity $\ge 2$.
\end{proposition}

\begin{proof} Assume that $s\ge 2$.
Denote the weight indicated in the proposition by
$$\mu :=2(s\omega_{i_1}+(s-1)\omega_{i_2}+\ldots +\omega_{i_s})-\epsilon_{i_1}+\epsilon_{i_s+1}$$
The bundle $\eta$ can be filtered so the associated graded $gr(\eta)$ is the direct  sum of line bundles, one for
each weight space of the $\fl$-module $M_{\mu}$. Their weights are as follows. 
Denote $\rho_P =\o_{i_1}+\ldots +\o_{i_s}$. Then the weights in $gr(\eta)$ are $\rho_P$ and $\rho_P -\alpha$ where $\alpha$ runs through the positive roots not in $\fp$.
We calculate the multiplicity of the weight $\mu$ in $S^2 (gr(\eta) )$. We have to count the cardinality of the set of pairs of weights in $gr(\eta)$ which add up to $\mu$. It is the set of sums of two roots $\beta_1 ,\beta_2$ (as to choose weights $\rho_P -\beta_1$, $\rho_P -\beta_2$), not in $\fp$ which add up to $\epsilon_{i_1}-\epsilon_{i_s+1}$, plus one (coming from the weights $\rho_P$ and $\rho_P -\epsilon_{i_1}+\epsilon_{i_s+1}$). This multiplicity is $i_s -i_1 +1$. Next we notice (using Bott-Borel-Weil) that the only occurrence of the weight $\mu$ in higher cohomology are the weights giving $\mu$ in
 $H^1 (S^2  gr(\eta) )$ and the multiplicity with which it occurs is the set of pairs of roots $\beta_1 =\epsilon_{i_1}-\epsilon_{j}$ and $\beta_2 =\epsilon_{j+1}+\epsilon_{i_s+1}$ such that $i_1<j\le i_s$ is not the last element in the corresponding interval
  $[i_{u}+1,\ldots ,i_{u+1}]$.
This multiplicity equals
$$(i_2-i_1 -1)+(i_3-i_2-1)+\ldots +(i_s -i_{s-1}-1)=i_s -i_1 -(s-1).$$
Thus the difference of multiplicities of $V_{\mu}$ in $H^0 (S^2 (gr(\eta)))$ and in $H^1(gr(\eta))$ is equal $s\ge 2$ which proves the Proposition.
\end{proof}
\end{proof}

\section{Irreducible cominuscule varieties}\label{irredchsssect}
Let $X=G/P_{\aio}\subset \BP V_{\oio}$ be an irreducible 
rank $r\geq 3$ cominuscule variety. Continuing the notation of \S\ref{isingsect}
with $\eta=(\hat TX)^*$,
$$
gr (\eta) = \cO_X(1)\oplus (\cO_X(1)\ot T^*X).
$$

Recall that   homogeneous vector bundles 
$E\ra G/P$   correspond to $\fp$-modules $M$, where
$E=G\times_PM$.  
In particular, the tangent bundle $TX$ corresponds   to the $\fp$-module $\fg/\fp$.

Recall further that irreducible $\fp$-modules are in one to
one correspondence with irreducible $\fl$-modules, which
are indexed by the set of  $\fl$-dominant weights,
which we denote $\Lambda^+_{\fl}$. 
Note that  $\l\in \Lambda^+_{\fl}$
iff $\l=a_1\o_1+\cdots + a_{\ell}\o_{\ell}$,
with $a_j\in \BZ$ and $a_j\geq 0$ for $j\neq i_0$,
where the $\o_j$ are the fundamental weights of $\fg$
and $P=P_{\aio}$. 

Let $M_{\l_1}$ denote $T^*_{[v_{\oio}]}X$ considered
as an $\ff:=[\fl,\fl]$-module.
A uniform, for all cominuscule varieties $X$, decomposition of $Sym(M_{\l_1})$ as an
$\ff$-module is given in  \cite{LMseries}.
These modules are exceptional in the sense of
Brion \cite{briont}, that is the symmetric algebra is
free. The formula is 
$$
S^d M_{\l_1}= \bigoplus_{j_1+2j_2+\cdots +rj_r=d}
M_{j_1\lambda_1 +\ldots +j_r\lambda_r}
$$
where $M_{\l_2}$ is the complement of $M_{2\l_1}$
in $S^2(M_{\l_1})$ and
 $M_{\l_j}=S^jM_{\l_1}\cap(S^{j-2}M_{\l_1}\ot M_{\l_2})
\subset M_{\l_1}^{\ot j}$.
In other words, consider the composition
$\d_j: S^{j-2}M_{\l_1} \ot S^2M_{\l_1}\ra
M_{\l_1}^{\ot j}\ra S^{j-1}M_{\l_1} \ot \La 2M_{\l_1}$,
then
$M_{\l_j}=\tker \d_j|_{S^{j-2}M_{\l_1}\ot M_{\l_2}}$.

\smallskip

\begin{remark}
The above correspondence among $\ff$-modules,
where generators of the symmetric algebra of
$M_{\l_1}$ under the Cartan  product correspond
to the prolongations of $M_{\l_2}$, 
 extends to the following correspondence
among $\fl$-modules:

{\it 
The generators of the ring of covariants
of $Sym(gr(\eta))$ correspond  to the
irreducible components of $gr(\xi)$.}

This correspondence is via the projective fundamental
forms. For any variety $Z\subset \BP V$ and
$z\in Z_{smooth}$, we have maps
$\BF_j: S^jT_zZ\ra N_zZ$. If we let 
$N_k=\tim \BF_j$, then 
$gr(\xi)=\oplus_k N_k^*$, and the 
$j$-th generator of 
the ring of covariants
of $Sym(gr(\eta))$ is
${}^t\BF_j(N^*_k)\subset S^jT_z^*Z$.
This clarifies the cryptic remarks on
p. 80 of \cite{LM0}.
\end{remark}

\smallskip

\begin{remark}
$M_{\l_j}$ admits the geometric interpretation of
the generators of the ideal of $\s_j(F/Q)\subset \BP M_{\l_1}^*$,
the variety of secant $\BP^{j-1}$'s to $F/Q=F[v_{\l_1}]$,
where $v_{\l_1}$ is a highest weight vector in $M_{\l_1}^*$.
See \cite{LMseries} for more information.
\end{remark}

\smallskip

Here is a table of the   rank $r$ cominuscule varieties,
  together with a description of the $F$-modules $M_{\l_j}$ and $T_{[v_{\oio}]}X$ as an $F$-module:

\begin{center}\begin{tabular}{||c|c|c|c|c||} \hline 
$X$  &  $G(k,n)$  & $G_{Lag}(n,2n)$ & $\SS_{2n}$ & $\QQ^n$\\ 
$G$  & $SL_n$  & $Sp_{2n}$ & $ Spin_{2n}$ & $SO_{n+2}$ \\  
$F$ &   $SL_k\times SL_{n-k}$  & $SL_n$ & $SL_n$ & $SO_n$\\
$T_{[v_{\oio}]}X$ & $M_{\o_{1}+\o_{n-1}}=(\cS^*\ot \cQ)_{[v_{\oio}]}$ &  $M_{2\o_1}=(S^2\cS^*)_{[v_{\oio}]}$ & 
$M_{\o_ {2}}=(\Lambda^2\cS^*)_{[v_{\oio}]}$ & $M_{\o_2}$\\ 
$M_{\lambda_j}$ & $M_{\o_{k-j}+\o_{k+j}}=(\La j \cS\ot \La j \cQ^*)_{[v_{\oio}]}$ & 
$M_{2\o_j}=(S_{{\scriptstyle 2\ldots 2}}\cS)_{[v_{\oio}]}$ & $M_{\o_ {n-2j}}=(\Lambda^{2j}\cS)_{[v_{\oio}]}$
& $\CC$ \\ 
$r$ & $\min (k,n-k)$ & $n$ & $\llcorner \frac{n}{2}\lrcorner$ &  $2$ \\ \hline   
\end{tabular}\end{center}

\medskip

\begin{center}\begin{tabular}{||c|c|c||} \hline 
 $X$  & $\BO\pp 2$ &     $G_{\o}(\BO^3,\BO^6)$       \\ 
  $G$  &$E_6$ &  $E_7$   \\ 
  $F$ &$ Spin_{10}$ &  $E_6$   \\
  $T_{[v_{\oio}]}X$  &$M_{\o_ 2}$ &  $M_{\o_6}$   \\
  $M_{\lambda_2}$ & $M_{\o_ 6}$ &   $M_{\o_1}$   \\
$M_{\lambda_3}$  & $0$  &   $\CC$  \\
$r$ & $2$&$3$   
\\ \hline \end{tabular}\end{center}

We caution the reader that the $\l$ in the
$M_{\l}$ are to be considered as highest
weights as $\ff$-modules, but the labeling
is as an element of the weight lattice of $\fg$.

Here $\cS,\cQ$ are respectively the
tautological subspace and quotient bundles.

\medskip 

We almost have a description of $Sym(gr(\eta))$ from this,
what is missing is the coefficient of the
weights on $\o_{i_0}$.
Because irreducible $\fl$-modules correspond to irreducible
$\ff$-modules equipped with an integer weight
on $\o_{i_0}$, 
adopting the notations that $M_{\mu}$ is the $\fl$-module with
highest weight $\mu$, 
and  $E_{\mu}$ is the corresponding irreducible
vector bundle on $X$, we may write
$$
S^d (T^*X\ot \cO_X(1))= \bigoplus_{j_1+2j_2+\cdots +rj_r=d}
E_{j_1\mu_1 +\ldots +j_r\mu_r}
$$
where $\mu_j=\l_j + m_j\o_{i_0}$ for
some integers $m_j$, which we now determine.

\begin{lemma} Notations as above.
$m_j=j-2$.\end{lemma}

\begin{proof}
There is a unique element $U_{i_0}$ of $\ft$,
called the {\it grading element} that
has the property $U_{i_0}(\a_j)=0$ if $j\neq i_0$
and $U_{i_0}(\aio)=1$. See e.g. \cite{yamaguchi}, \S 3.1.
In particular,     $TX$  is an eigenspace
for the action of $U_{i_0}$ with eigenvalue one,
and $S^j(T^*X)$ is an eigenspace with
eigenvalue $-j$ and $S^j(T^*X\ot \cO(1))$
is an eigenspace with eigenvalue
$-j+j(c\inv)_{i_0,i_0} $
where $c\inv$ is the inverse of the Cartan matrix.
Thus
$$
-j+j(c\inv)_{i_0,i_0} 
=U_{i_0}(\mu_j)=U_{i_0}(\l_j)+m_j(c\inv)_{i_0,i_0}
$$

The lemma thus reduces to showing
\begin{equation}\label{Zid}
U_{i_0}(\l_j)=
2(c\inv)_{i_0,i_0} -j
\end{equation}
which can easily be checked on a case by case basis.
\end{proof}

\begin{remark}
A uniform and conceptual proof 
of \eqref{Zid} is possible, but it would take us too far afield here.
In particular, note that $\mu_1=\oio -\aio$.
\end{remark}

\begin{example} To verify \eqref{Zid} in the
case of $D_n/P_n\subset \BP V_{\o_n}$,  we have $i_0=n$, $\l_j=\o_{n-2j}$,
$U_{n}(\o_i)=(c\inv)_{i,n}=\frac i2$ for $i< n-1$
and $U_n(\o_n)=\frac n4$. We indeed
have $\frac{n-2j}2=2(\frac n4) -j$.
\end{example}

\begin{example} To verify \eqref{Zid} in the
case of $A_n/P_k\subset \BP V_{\o_k}$,  we have $i_0=k$, $\l_j=\o_{k-j}+\o_{k+j}$,
$U_{k}(\o_{k-j})=\frac{(k-j)(n-k+1)}{n+1}$,
$U_k(\o_{k+j})=\frac{k(n-(k+j)+1)}{n+1}$,
$(c\inv)_{k,k}=\frac{k(n-k+1)}{n+1}$  
and we verify
$$
\frac{(k-j)(n-k+1)}{n+1} +\frac{k(n-(k+j)+1)}{n+1}
=2\frac{k(n-k+1)}{n+1} -j.
$$
\end{example}

\smallskip

Using that $\mu_j=\l_j+(j-2)\oio$ and
$S^d(gr(\eta))=\oplus_{k=0}^dS^k(T^*X\ot \cO(1))\ot \cO(d-k)$,
we conclude:

\begin{proposition}Notations as above.
$$
S^d(gr(\eta))=  \bigoplus_{a_1+2a_2+\cdots +ra_r  \leq d}  
E_{a_1\mu_1 +\ldots +a_r\mu_r+(d-\sum_{j=1}^rja_j)\oio}.
$$
\end{proposition}

\medskip

Recall that the Bott-Borel-Weil theorem
implies that for irreducible homogenous vector bundles
$E_{\mu}\ra G/P$, 
\begin{enumerate}
\item $\oplus_kH^k(E_{\mu})$ is
an irreducible $G$-module (and in particular
is nonzero in at most one degree $k$), 
\item writing
$\mu=\sum_ia_i\o_i$, if all $a_i\geq 0$
(i.e. if $\mu\in \Lambda^+_{\fg}$), then
$H^0(E_{\mu})=V^*_{\mu}$, 
\item  if all
$a_i$ but $a_{i_0}$ are non-negative,
and $a_{i_0}=-1$ then there is no cohomology,
\item  if all
$a_i$ but $a_{i_0}$ are non-negative,
  $a_{i_0}<-1$, and moreover
$\s_{\aio}.\mu\in  \Lambda^+_{\fg}$ then $H^1(E_{\mu})=
V^*_{\s_{\aio}.\mu}$, where
$V_{\nu}$ is the $\fg$-module
with highest weight $\nu$, 
(Here $\s_{\aio}.\mu$ denotes the
affine action of the Weyl group,
$\s_{\aio}.\mu= \s_{\aio}(\mu +\rho)-\rho$,
where $\s_{\aio}$ is reflection in the
hyperplane orthogonal to $\aio$.)
\end{enumerate}

\begin{remark} One obtains the dual modules as cohomology
groups above because of our convention of deleting
negative root spaces to define our parabolic subalgebras.
\end{remark}

\medskip

In our case, for each $E_{\mu}$ that
appears, the only possible negative
coefficient is that of $\oio$.
Moreover,
for $j>1$, $\s_{\oio}.(\l_j)=\l_j$ for $j>1$
because   all the $\l_j$
except for $\l_1$ are
orthogonal to $\oio$. (This  can be verified
case by case, but it is also
 a consequence of what is often called \lq\lq Kostant's cascade\rq\rq .)  
For $j=1$, we have 
$\s_{\aio}.(a_1\l_1+c\oio)=(a_1+c+1)\l_1-(2+c)\oio$.
This last assertion follows immediately from
the observation that $\mu_1=\oio-\aio$.
In particular, if $E_{\mu}$ appearing
in $S^d(gr(\eta))$ is neither ample, nor
has no cohomology, then $\s_{\aio}.\mu\in \Lambda^+_{\fg}$.
In summary:

\begin{proposition}\label{Ecohprop} 
Recall the notations that $M_{\mu}$ is the $\fl$-module with
highest weight $\mu$, 
  $E_{\mu}$ is the corresponding irreducible
vector bundle on $X=G/P_{i_0}\subset \BP V_{\oio}$,     and $V_{\nu}$ is the $\fg$-module
with highest weight $\nu$.  
If $a_1+2a_2+\cdots + ra_r\leq d$, then,
letting 
\begin{align*}
\mu&=a_1\lambda_1 +\ldots +a_r\lambda_r+(d-2\sum_{j=1}^r a_j)\oio\\
&= a_1\mu_1 +\ldots +a_r\mu_r+(d- \sum_{j=1}^r ja_j)\oio,
\end{align*} we have
 
\begin{enumerate}
\item $E_{\mu}$
is ample with $H^0(E_{\mu})=V^*_{\mu}$
when $d-2\sum_{j=1}^r a_j\geq 0$,

\item $E_{\mu}$ has
no cohomology when $d-2\sum_{j=1}^r a_j=-1$,

\item  $E_{\mu}$ has 
$
H^1(E_{\mu})
=V^*_{\s_{\oio}.\mu} 
$
when $d-2\sum_{j=1}^r a_j<-1$. 
\end{enumerate} 
\end{proposition}

Note that
$$\s_{\oio}.\mu=
(a_1+ d-\sum_{j=1}^r2a_j +1)\lambda_1+ a_2\l_2 +\ldots +a_r\lambda_r-(2+d-2\sum_{j=1}^r a_j)\oio.
$$

\begin{remark} Using the $\mu_j$ and $\oio$ has
the advantage that these are the actual highest weights
of the primitive $\fl$-modules that show up in the
decomposition of $Sym(gr(\eta))$, while using the $\l_j$ and $\oio$
has the advantage that
all but
$\l_1$ are orthogonal to $\oio$, and
 (except for the
case of $\fa_n$) all the $\l_j$  
are orthogonal to each other, and the $\l_j$   are
  fundamental weights for $\fg$, with the
exception of $\fg=\fa_n$ where they are
sums of fundamental weights.
\end{remark}

Now that we have determined the cohomology of
$Sym(gr(\eta))$ we turn to $Sym(\eta)$.
At most the cohomology groups appearing
in $Sym(gr(\eta))$ can appear, but there can
be cancellation. Note first that for a given
$\mu$, $E_{\mu}$ appears at most once in $Sym(gr(\eta))$.
Moreover, for the $E_{\mu}$ appearing in $S^d(gr(\eta))$ with 
$H^1(E_{\mu})$ nonzero, the bundle $E_{\mu'}$ with 
\begin{align*}
\mu'&=(a_1+(d-2\sum_{j=1}^r a_j)+1)\lambda_1+ a_2\l_2 +\ldots +a_r\lambda_r+(-2-d+2\sum_{j=1}^r a_j)\oio\\
&
=(a_1+(d-2\sum_{j=1}^r a_j)+1)\mu_1+ a_2\l_2 +\ldots +a_r\mu_r+(-1-2d-\sum_{j=1}^r(j-2) a_j)\oio
\end{align*}
also
appears in $S^{d}(gr(\eta))$ and of course
$H^0(E_{\mu'})=H^1(E_{\mu})$.
It remains to show that these terms cancel when
one passes to $H^*(Sym (\eta))$.

In order to prove that the matching terms cancel out in the spectral sequence we use the technique of \cite{Orub}.
The essential point is that 
$$Ext^1(\cO(1),T^*X\ot \cO(1))=H^1(Hom(\cO(1),T^*X\ot \cO(1)))
$$
(see \cite{Hart}, proposition 6.5),
and $Hom(\cO(1),T^*X\ot \cO(1))\simeq T^*X$.
Now $T^*X= E_{\l_1-2\oio}$ and applying the
Bott-Borel-Weil theorem again, we see
$H^1(E_{\l_1-2\oio})=\BC$. Thus there is a unique
(up to scale)  nontrivial extension.

The quiver representation of a quiver $\cQ_X$ defined in \cite{Orub} corresponding to $S^d  \eta $ has one dimensional spaces attached to vertices with the highest weights of $L$-modules $  M_{(d-\sum ja_j)\oio +a_1\mu_1 +\ldots +a_r\mu_r}$ for $a_1+2a_2+\cdots +ra_r  \leq d$. The arrows connect the weight of 
$M_{\mu}$ to that of $M_{\mu'}$ where $\mu,\mu'$ are
as above.
Then,
noting that
$\eta$ is indeed a nontrivial extension because
it is acted on nontrivially by $\fn$,  \cite{Orub},  Proposition 6.7 assures that the connecting homomorphism between two
cancelling terms is nonzero. 

Recall that for an algebra $\cA$ that has
the structure of a $\fg$-module,
  the {\it ring of covariants} of $\cA$ is the
set of elements of $\cA$ annihilated by all positive
root vectors in $\fg$ (or, if working with
a corresponding algebraic group $G$, 
the elements invariant under the action of
a unipotent radical of $G$). Another perspective
is that the ring of covariants is the generators
of $\cA$  as an algebra with the multiplication
by Cartan product instead of its usual multiplication.

\begin{theorem}\label{irrcoordgens} Let $X=G/P_{i_0}
 \subset \BP V_{\oio}$ be   rank
$r\geq 3$ cominuscule variety.    Let
$K[\t (X)]$ denote the homogeneous coordinate ring of $\t (X)$. Then,
 continuing the notation of above

\begin{enumerate}
\item
 
$$
K[\t(X)]_d = \bigoplus_{  2\sum_{j=1}^r a_j\leq \tmin\{d,\sum_{j=1}^r ja_j\} }
 V^*_{(d- \sum_jja_j)\oio +a_1\mu_1 +\ldots +a_r\mu_r}  
$$

\item  
The ring of covariants of
$K[\t(X)]$    is generated by the modules
$V_{\omega_{i_0}}^*$, $V^*_{i\mu_1 +\mu_s}$ with $1\le i\le s-2$, and $3\le s\le r$,
and $V^*_{\mu_2},\ldots ,V^*_{\mu_r}$. Thus the ring is generated in degrees $\le 2(r-1)$.
  \end{enumerate}
\end{theorem}

\begin{remark} The marked Dynkin
diagram describing the module $V_{\mu}^*$ is  
the marked Dynkin diagram of $V_{\mu}$ reflected
by the $\BZ_2$-symmetry of the diagram.
\end{remark}

\begin{proof}To prove 
the first assertion,
  by Theorem \ref{weymanthm}.2
and the preceeding paragraph,
we just need to calculate $H^0(S^d(\eta))$.
By  Proposition \ref{Ecohprop}, 
a module  
$V^*_{(d- \sum_jja_j)\oio +a_1\mu_1 +\ldots +a_r\mu_r}$
is in
$H^0(S^d(gr(\eta))$ if $2\sum a_j\leq d$, and
by the discussion above, to
obtain $H^0(S^d(\eta))$ we must
subtract the modules $V^*_{\mu}$
such that $(\s_{\aio}.)\inv \mu$
also occurs in $S^d(gr(\eta))$.
Since $(\s_{\aio}.)\inv (a\l_1+c\oio)
=(a+c+1)\l_1+(-2-c)\oio$
we need to subtract the modules with
$$
(a_1+(d-2\sum_ja_j)+1)+2a_2+\cdots +ra_r
\leq d,
$$  i.e.,
we require
$(a_1 -2\sum_ja_j +1)+2a_2+\cdots +ra_r
> 0$,   i.e.,
$2\sum a_j\leq \sum ja_j$.

\smallskip

For the second
assertion,  it is clear that   $V^*_{\omega_{i_0}}$ is among the generators. 
To prove that the other generators are as described in the statement, 
  we use
induction on $s:=\tmin \{ i\mid i>1, a_i>0\}$.  
Consider an $r$-tuple $a=(a_1\hd a_r)$ such that
$a_2\hd a_{s-1}=0$ satisfying $2\sum a_j\leq ja_j$. Set $k=\tmin (a_s,\llcorner \frac{a_1}{s-2}\lrcorner)$
and subtract $k(s-2,0\hd 0,1,0\hd)$. Either one obtains $a_s=0$
and we may go to the next step of the induction or
$a_s>0$ and $a_1<s-2$. But such a vector $a$ is a non-negative
linear combination of vectors $(a_1,0\hd 0,1,0\hd )$ and the vectors
in the basis for larger $s$.
\end{proof}

To determine generators of the ideal we must calculate the modules
$H^d(\La{d+1}\xi)$. Here $\xi= N^*_X(1)$.
Unfortunately even decomposing $\La i M_{\lambda_j}$ in general is a difficult
problem, which is why we are only able to determine explicit generating
modules in a few special cases.

\section{Multiple embeddings}\label{multsect}

In this section we generalize the results of \S\ref{irredchsssect} to multiple embeddings of
$X$. The results become   easier, because all higher cohomology of $Sym (gr(\eta))$ vanishes.
Assume that $X$ is embedded into $V_{N\o_{i_0}}^*$ by the $N$-tuple embedding with $N\ge 2$.

Using  the notation of \S\ref{isingsect},  
with $\eta=(\hat TX)^*$,
$$
gr (\eta) = \cO_X(N)\oplus (\cO_X(N)\ot T^*X).
$$
This implies
\begin{proposition}Notations as above.
$$
S^d(gr(\eta))=  \bigoplus_{a_1+2a_2+\cdots +ra_r  \leq d}  
E_{a_1\mu_1 +\ldots +a_r\mu_r+(Nd-\sum_{j=1}^rja_j)\oio}.
$$
\end{proposition}
By repeating the reasoning from the previous section we get

\begin{theorem}\label{irrmultgens} Let 
$X=G/P_{i_0}\subset \BP V_{N\o_{i_0}}$ be a   irreducible rank
$r\geq 3$ generalized cominuscule variety, embedded by an $N$-tuple embedding with $N\ge 2$.    Let
$K[\t (X)]$ denote the homogeneous coordinate ring of $\t (X)$. Then,
 continuing the notation of above

\begin{enumerate}
\item
 
$$
K[\t(X)]_d = \bigoplus_{  \sum_j ja_j\le d} 
 V^*_{(Nd- \sum_jja_j)\oio +a_1\mu_1 +\ldots +a_r\mu_r}  
$$

\item  
The ring of covariants of
$K[\t(X)]$    is generated by the modules
$V^*_{N\omega_{i_0}}$, 
and $V^*_{(N-1)\o_{i_0}+\mu_j}$, for $1\le j\le r$. Thus the ring is generated in degree 1.
However, since in degree 1 we have more than one representation, the embedding 
of $\t(X)$ into $V^*_{\o_{i_0}}$ is not
linearly normal.
  \end{enumerate}
\end{theorem}
\begin{proof}
The result follows at once by observing that all bundles $ E_{a_1\mu_1 +\ldots +a_r\mu_r+(Nd-\sum_{j=1}^rja_j)\oio}$ are now ample.
\end{proof}

\begin{corollary}
Let $X\subset \BP V$ be a rank $\geq 3$ cominuscule variety.
Then $\t(X)$ is not quadratically normal.
\end{corollary}

\section{Reducible cominuscule varieties}\label{reduciblecominsect}
Let $X=Seg(X_1\times \cdots \times X_m)\subset
\BP V=\BP (W_{\o_{i_0}^1}\ot \cdots \ot W_{\o_{i_0}^m})$
 where each $X_i=G^i/P_{ \a_{i_0}^i}\subset \BP W_{\o_{i_0}^i}$
is a   rank $r_i$ 
cominuscule variety so the rank of $X$ is $r:=r_1+\cdots +r_m$.
(We leave the case of non-minimally embedded factors to the reader.)
Write $\cO_X(p_1\hd p_m)=\cO_{X_1}(p_1)\ot \cdots \ot \cO_{X_m}(p_m)$
where we have omitted the pullback maps from the notation.
Then
$$
gr(\hat T X)= 
\cO_X(-1\hd -1)\op \op_{j=1}^m \cO_X(-1\hd -1)\ot TX_j.
$$
 
We adopt the notation ${}^sM_{\lambda^s_j}=M_{\lambda^s_j}(X_s)$ 
and $\eta_s=\eta (X_s)$ following the notation of \S\ref{irredchsssect} for the irreducible cases, in particular
$a^{s}_{j_{s}}$ corresponds to $a_{j_{s}}$ of
the $s$-th factor and similarly for
$\o^s_{i_0}$, $\mu^s_j$ etc...
Then
$$
S^d(gr(\eta))=\bigoplus_{p_1+\cdots +p_m\leq d}
\cO_X(d-p_1\hd d-p_m)\otimes 
S^{p_1}(gr(\eta_1))\ot \cdots \ot S^{p_m}(gr(\eta_m)).
$$

\begin{lemma}\label{41} Notations as above.

\begin{enumerate}

\item $H^k(S^d(gr(\eta)) =0$ for $k>1$. 

\item The modules of the form $\ot_{s }  {}^sM_{\mu^s}$  are in the positive
cone of the Weyl chamber are those whose $s$-th component  
is in the positive cone for each $s$.

\item The modules  appearing
in $S^d(gr(\eta))$ with no $H^0$ term are those where
$d<2\sum_{j=1}^{r_s}a^s_j$ for some $s\in \{ 1\hd m\}$.
Note that this can occur for at most one such $s$ and that such
a term will contribute a module to $H^1$.

\item The modules appearing in   $H^1(gr(Sym(\eta)))$ all appear
in  $H^0(gr(Sym(\eta)))$ with the same
multiplicity. These terms cancel when one passes to
$H^*(Sym(\eta))$.
\end{enumerate}
\end{lemma}

Lemma \ref{41} follows from proposition
\ref{Ecohprop}  by observing that the
Weyl group acts independently on each  factor and at most one
factor can fail to be ample.

\begin{theorem}\label{redcoordgens} Let $X=Seg(X_1\times \cdots \times X_m)\subset \BP V=\BP (W_{\oio^1}\ot \cdots \ot W_{\oio^m})$ be a homogeneously embedded   rank
$r=r_1+\cdots + r_m$ cominuscule variety with $r\geq 3$.  Let
$K[\t (X)]$ denote the homogeneous coordinate ring of $\t (X)$.
  Then

\begin{enumerate}
\item Let $p_s=a^s_1 +2a^s_2+\ldots +r_s a^s_{r_s}$.

$$
K[\t(X)]_d = \bigoplus'
\bigotimes_{s=1}^m
  {}^sV^*_{(d-\sum_{j=1}^{r_s}ja^s_j)\o^s_{i_0}+a^s_1\mu^s_1
+\cdots + a^s_{r_s})}.$$
 
The sum $\oplus '$ is over sets $(a^s_j)$
such that 
$\forall s,\ 2\sum_{j=1}^{r_s}a^s_j
\leq\tmin\{ d,\sum_{j=1}^{r_s}ja^s_j+p_1+\cdots +\hat p_s+\cdots
+ p_m\}$ and 
$p_1+\cdots +p_m\leq d$.

\item  The generators, which  (aside from
those of type (i.))   we label by the sets of
integers $a^s_j$     come in four  types:

(i.) Fix $s$ and for each $e\ge 0$ take the generators in degree $e$ of $K[\t (X_s )]$  on the $s$-th coordinate 
(listed in Theorem \ref{irrcoordgens})
tensored with the representation ${}^tV^*_{e\o_{i_0}}$ on the remaining coordinates,

(ii.) Fix $s_1,s_2\in \{ 1\hd m\}$ and $j_2\in \{1\hd r_{s_2}\}$, we have $j_2\geq a^{s_1}_1>0$,
$a^{s_1}_{\rho}=0$ for $\rho>1$,
$a^t_j=0$ for all $t\neq s_1$
except for $a^{s_2}_{j_2}=1$.

(iii.)Fix $s_1,s_2$, $a^{s_1}_1=a^{s_2}_1=1$ and all other
$a^t_j$ are zero.

(iv.) Fix $s_1,s_2,s_3$, $a^{s_1}_1=a^{s_2}_1=a^{s_3}_1=1$ and all other
$a^t_j$ are zero.

Thus  $K[\t(X)]$ is generated in degrees up
to $\tmax\{3, 2r_0\}$ where $r_0=max_sr_s$.

\end{enumerate}
\end{theorem}

The proof is similar to the irreducible case.

\begin{example}[Segre varieties]\label{segex} Consider $m$ factors of $\BP^1$,
$Seg(\pp 1\times \cdots\times \pp 1)\subset
\BP (K^2\otc K^2)$. The representations occurring in
$K[\t(X)]_d$ are all modules that are tensor products of Schur functors with
partitions of $d$ of length at most two.
For each $s=1,\ldots ,m$ the weight $ \oio^s$ 
corresponds to the partition $(1,0)$ and the weight $\lambda^s_1$
corresponds to $(0,1)$.
The generators of $K[\t(X)]$ are as follows. There is a representation with the  weight
$\otimes_{s=1}^m \lambda^s_0 =\otimes_{s=1}^m (1,0)$ in degree one, and the representations with the weights
$\n_I$ for any $I\subset [1\hd n]$ $|I|=2$ or $3$,
where for $|I|=i$ we have $(\n_I)_s=(i-1,1)$ for $s\in I$ and $(\n_I)_s=(i,0)$ for $s\notin I$. Since the rank of each $X_i$ is one, we have that all generators of $K[\t(X)]$ occur in degree at most three. It will follow from the results of \S\ref{degensect}  that
$I(\t (X))$ is generated in degrees at most six.
\end{example}

\section{The degeneration argument}\label{degensect}

We use the notation of \S 15 of \cite{Gr}.
Since in \cite{Gr} algebraic groups are used instead
of Lie algebras when discussing weights etc..., for
this section only we use $\Lambda_G$ instead of $\Lambda_{\fg}$,
  although
since we are in characteristic zero, we could have
just as well used Lie algebras.  Let $G$ be a linearly reductive group, $T$ a maximal torus and $U$ a unipotent radical.
For an algebra $A$ with rational $G$-action Grosshans (Lemma 15.1) constructs a homomorphism $h:\Lambda_{G}\rightarrow {\bf Z}$ satisfying the properties  

\begin{itemize}
\item{a)} $h(\omega )\in \BZ_{\geq 0}$ when $\omega\in \Lambda^+_G$,
\item{b)} if $\chi >\chi^\prime $(i.e. the difference is a sum of positive roots), then $h(\chi )>h(\chi^\prime)$,
\item{c)} $h(g_j\chi )=h(\chi )$ where $\lbrace g_j\rbrace$ is the set of representatives of cosets of $G$ with respect to
the connected component of the identity $G^0$ (this is trivial for connected $G$),
\end{itemize}

For an algebra $A$ with a rational $G$-action Grosshans defines
$$A_n=\lbrace a\in A\ |\ h(\chi )\le n {\rm \ for\ all\ weights\ }\chi{\rm\ of\ }T {\rm \ in\ the\ span\ }\langle G.a\rangle\rbrace .$$

We define 
$$gr(A):= \oplus_{n\ge 0} (A_n /A_{n-1}).$$

This is a commutative algebra with a rational $G$-action, with the product induced by the product in $A$. The algebras $A$ and $gr(A)$ have the same algebras of $U$-invariants.
Define  
$$D:=\sum_{n\ge 0} A_n x^n\subset A[x].$$

The algebra $D$ has a rational $G$-action and it has the following properties:
\begin{itemize}
\item{d)} $D/xD =gr(A)$,
\item{e)} $D[{\frac 1x}]=A[x, {\frac 1x}]$.
\end{itemize}

Theorem 15.14 in \cite{Gr} implies:

\begin{theorem}  Let $i:K[x]\rightarrow D$ be an inclusion. Then $D$ is flat over $K[x]$. The fiber of $i$ over a maximal ideal $(x-\alpha)$, $\alpha\in K^*$, is isomorphic to $A$, and the fiber over $(x)$ is isomorphic to $gr(A)$.
\end{theorem}

Another way to think about $gr(A)$ is that it is $A$ with the product deformed so only the Cartan piece of the product of $A$ is retained.

We apply the theorem to multiplicity free algebras.

\begin{theorem}\label{degreethm}  Assume that $A$ is multiplicity free graded domain with rational $G$-action preserving the grading. Assume $gr(A)$ is generated by representations of degree $\le d$.
Then the defining ideal of $A$ is generated in degrees $\le 2d$.
\end{theorem}

\begin{proof} Let
$\Theta:=\{ \l\in \Lambda_G^+\mid V_{\l}\subset A\}$.
Note that $A$ is multiplicity free and a domain, so
$V_{\l},V_{\mu}\subset A$ implies $V_{\l+\mu}\subset A$ thus
$\Theta$ is an abelian sub-semi-group of $\Lambda_G^+$.
Consider $gr(A)$. This is an algebra that additively is
$$gr(A)=\oplus_{\lambda\in\Theta}V_\lambda$$
with the product given by Cartan product. By the previous theorem $gr(A)$ is a special fiber  of a flat deformation with general
fiber $A$. Introduce a new degree on $gr(A)$ by setting the degrees of the generators of $\Theta$ to one. Now by an
unpublished theorem of Kostant  \cite{Gar}  $gr(A)$ has relations in degrees $\le 2$. This means the original degrees of these relations are $\le 2d$. But there is a  presentation of the general
fiber  given by the  generators and relations in the same degrees as that of a special
fiber. This proves our statement.
\end{proof}

\section{Further information on the ideals}\label{exsect}\label{examsect}

First note that $\t(X)\subseteq \s (X)$ and, as discussed in
\cite{LWsecseg}, in most cases there is a {\it subspace variety}
or {\it rank variety} containing $\s(X)$. Thus we may study the
equations of $\t (X)$ by first understanding certain \lq\lq primitive\rq\rq\
 cases
  and then the ideals of the rank varieties themselves.
The ideals of rank varieties are possible to determine by the   method of \cite{weyman} in many cases.
Assuming we have both a set of generators of the ideal of a primitive
case and of the relevant rank varieties, one must still determine which
generators of the ideal of the rank variety become redundant
when considered as members of the ideal of $\t (X)$.

\subsection{Grassmannians} 
Consider $G(r,N)$. Write $V^*=\BC^N$.
In this case the primitive varieties are $G(r,2r)$ and the
relevant
rank variety in $\La r V^*$ is   
$$Z_{r,N}=\{T\in \La r V^*\mid \exists W\subset V^*, \tdim W=2r,\ 
T\in \La r W\}.$$ 
These rank varieties are
discussed in \S 7.3 of \cite{weyman}.

\smallskip

  Consider $G(3,N)$. Here the primitive case is
$G(3,6)$ and 
$\t(G(3,6))$ is a quartic hypersurface whose equation is the
unique occurance of the trivial representation of $\fsl_6$ in
$S^4(\La 3\BC^6)$ corresponding to the partition $S_{222222}(\BC^6 )=S_{2^6}(\BC^6 )$.
Thus the ideal of $\t (G(3,V^*))$ is spanned by  $S_{2^6}V$
and the generators of the ideal of the subspace varieties.
For $G(3,7)$ the ideal of the variety
of tensors of rank $6$ is generated
by $S_{3111111}V$. For $G(3,8)$ the variety of tensors of rank $\le 7$
is generated by $S_{
3,2^2,1^5} V$ but  in \cite{LWsecchss} we show that this is in the ideal
generated by $S_{2^6}V$ and $S_{3,1^6}V$.
(See \cite{LWsecchss} for proofs of the assertions regarding
the generators of these rank varieties.)

\begin{theorem} The ideal of $\t(G(3,V^*))\subset \BP( \La 3 V^*)$
is generated by $S_{2^6}V$ in degree four and
$
S_{3,1^6}V$ in
degree three.
\end{theorem}

For $G(4,8)$ we calculated the Euler characteristics
of the bundles    $ \La{d+1}\xi$ directly.
In small degrees we can recover
$H^d(\La{d+1}\xi)$ from the Euler characteristic
to prove:
  
\begin{theorem}
The ideal of $\t(G(4,8))$ has, among its generators,
 $S_{1^8}V$ in degree two, $S_{2^6}V$ and $S_{3^2,1^6}V$
in degree three and $S_{4,2^6}V$ in degree four.
\end{theorem}

We expect that these modules in fact generate the ideal.

For $r>4$ the calculation becomes more difficult.

\subsection{Legendrian varieties}
These are the cases where
$X\subset \BP V$ is
$v_3(\pp 1)$,
$Seg(\pp 1\times \pp 1\times \pp 1)$,
$G_{Lag}(3,6)=C_3/P_3$, $G(3,6)$, $\BS_6=D_6/P_6$,
$E_7/P_7$, and $Seg(\pp 1\times Q)$ where $Q$ is a quadric
hypersurface. Here $\t(X)$ is a quartic 
hypersurface. The ring of covariants is free and has
a uniform description: it is generated by
$V$ in degree one, $\fg$ in degree two, $V$ in degree three
and $V_2$ in degree four, where 
$V_2\subset \La 2 V$ is the complement of the line spanned
by the symplectic form. See \cite{LMseries} for more
details.

\subsection{Scorza varieties}\label{scorza}

The homogeneous varieties $X$ with $\t (X)$ degenerate (i.e.,
of   dimension less than $2\tdim X$) coincide exactly
with the rank 2  cominiscule varieties:
Segre varieties,
$Seg(\pp a\times \pp b)\subset \BP(K^{a+1}\ot K^{b+1})$,
Grassmannians of two-planes
$G(2,n)\subset \BP (\La 2 K^n)$, the Cayley plane $E_6/P_6\subset \BP^{26}$,
quadric hypersurfaces and
the spinor variety $D_5/P_5$ (a component of 
the Grassmannian of $5$-planes in $K^{10}$
  isotropic for a quadratic form).
The 
ideals of the tangential
varieties in these cases are either
empty (when $\t (X)=\BP V$) or are generated by cubics.

In all these cases the secant variety of $X$  coincides with the tangential variety.
These cases also have a uniform description as rank one
elements in Jordan algebras with the tangential varieties
as the rank at most two elements.
The ideal of $\t(X)$ is generated in degree three by the
three by three minors in the Jordan algebra, the algebra
of covariants of $K[\t(X)]$ is free and generated
by the minors of the Jordan algebra of various sizes
(other than three by three). See \cite{LMseries} for
more details.

\subsection{Spinor varieties}
For $D_4,D_5$ the tangential varieties $\t (D_n/P_n)$ are
the ambient spaces, and these modules are exceptional in
the sense of Brion \cite{briont}, in that the ring of covariants is free.
The generators are respectively $V_{\o_3}$ in degree one
and the trivial representation in degree two for $n=4$ and
$V_{\o_4}$ in degree one and $V_{\o_1}$ in degree two for
$n=5$.

For $n=6$ we are in the Legendrian case.

For $n=7$,
checking the generators of the lattice of weights occurring in $K[\t(X)]$ we see that   the ring  of covariants is still free.

For   $n\geq 8$ it is no longer free. In fact for $n=8$
among the generators of $K[\t(X)]$ are
$V_{\o_2+\o_7}$ in degree three,
$V_{\o_6+\o_7}$ in
degree five   $V_{\o_2+\o_4}$ in degree three and
$V_{2\o_7}$ in degree four. The Cartan products of the
first two and last two representations are the same.

\subsection{Segre varieties}\label{segsect}

We collect all results about the tangential varieties of the Segre embeddings of products of projective spaces.
We use the geometric technique and the deformation argument from \S\ref{degensect}.
We start with the case of products of projective lines.

\begin{theorem} Let $X=Seg(\BP A_1^*\times \cdots \times \BP A_m^*)
\subset \BP(A_1^*\ot \cdots \ot  A_m^*)$ where $A_j$ is a vector space of dimension $2$ for $j=1,\ldots ,m$. Then

\begin{enumerate}

\item  The ideal of $\t (X)$ is generated in degree at most $6$.

\item The last term in the minimal free resolution of
$K[\t(X)]$ is
$\otimes_{j=1}^m S_{(2^{m-1}-2,2^{m-1}-m+1)}A_j$.
\end{enumerate}
\end{theorem}

\begin{proof}
The first statement follows from Example \ref{segex} and from Theorem \ref{degreethm}. The second statement follows from
Theorem 5.1.2 in \cite{weyman}. The top term in the resolution is easily seen to come from the cohomology of top exterior power of $\xi$. Then the result follows from Serre's theorem on cohomology of line bundles on   projective space.
\end{proof}

Now we pass to the general case where
 $X=Seg(\BP A_1^*\times \cdots \times \BP A_m^*)
\subset \BP(A_1^*\ot \cdots \ot  A_m^*)$ with $\tdim A_j=a_j+1$ for $j=1,\ldots ,r$.

\begin{theorem}\label{73} For Segre varieties $X=Seg(\pp {a_1}\times \cdots \times \pp{a_m})$, the tangential variety  $\t(X)$ is
arithmetically Cohen-Macaulay.
\end{theorem}

\begin{proof}
In the following proof we   use a relative
version of the machinery in \cite{weyman}. We will be  
terse here because
 a very similar argument with more detail   is in \cite{LWsecseg},
\S 5.

Let $Sub_{r\hd r}$ denote the rank or subspace variety whose desingularization
is given by the rank $2^m$ tautological subspace bundle   
$\cR_1\otc \cR_m\ra \Pi_{j=1}^m G(2, A_j^* )=B$ (see \cite{LWsecseg}, \S 3).
Note that $\tdim Sub_{r\hd r}=2^m +\sum_{j=1}^m 2(a_j -1)$. We have
$\t (X)\subset Sub_{r\hd r}$ (as we even have the secant
variety of $X$ contained in $Y$, see e.g., \cite{LWsecseg}).

The variety $Sub_{r\hd r}$ has a desingularization that allows one to apply the geometric technique from \cite{weyman}. In the notation of Theorem \ref{weymanthm},  we take  
$\eta =\cR_1^*\otc \cR_m^*$ and $\xi:= Ker ((A_1\otimes\ldots\otimes A_m\ot \cO_B)\rightarrow \eta )$.

 We consider the sheaf of algebras $\cB:=  Sym (\eta )$.  
We show the hypotheses of Theorem
\ref{weymanthm} are satisfied. We need   the following lemma from \cite{LWsecseg}:

\begin{lemma}\label{nohighercoh}[\cite{LWsecseg}, Lemma 5.2] 
Let $\pi_j=(p_{j,1}\hd p_{j,r})$
 be partitions.  Consider the sheaf
$${\cM}:= \otimes_{j=1}^m S_{\pi_j }{\cR}_j^*\otimes {\cB}.$$
  
\begin{enumerate}
\item Assume that $p_{j,1}\ge -a_j +1$ for $1\le j\le m$. Then     ${\cM} $ is acyclic.
\item Assume that $ p_{j,1}\ge 0$ and $p_{j,1}\le r^{m-1}-r$ for $1\le j\le m$. Then the $Sym(A_1\otc A_m)$-module  $H^0(B,{\cM})$, which is
supported in $Sub_{r\hd r}$,
is a maximal Cohen-Macaulay module.
\end{enumerate}
\end{lemma}

Now we use the desingularization of   ${\hat\t}(X)$ by
$\hat TX$. It is a vector bundle of rank $m+1$ which is a factor of the bundle $\cR_1\otc \cR_m$ defining the desingularization of $Sub_{r\hd r}$. To estimate higher direct images we first   analyze the finite free resolution of   $\hat \t(Seg( {\BP}^1
\ctimes \BP^1))$
($m$ copies) in the relative setting (taking ${\cR}_j$ instead $A_j^*$).

We apply the above lemma to this resolution.
By Theorem \ref{73},(2) each term in this resolution satisfies conditions  (1) and (2) of   Lemma \ref{nohighercoh}. Thus each term has no higher direct images and its sections form a maximal Cohen-Macaulay module supported in $Y$. This proves the vanishing of higher direct images of the structure sheaf of $Z$ and thus proves the  rational singularities.
 Taking resolution of each term in that complex and using an iterated mapping cone construction we get a nonminimal resolution of the coordinate ring of tangential variety whose length equal its codimension. This resolution  implies that the coordinate ring $K[\t(X)]$ is Cohen-Macaulay.
\end{proof}

Calculating 
the Euler characteristic of $ \Lambda^{d+1} gr(\xi) $
in low degrees we uncover certain generators of
$I(\t (X))$ which we expect to generate the ideal.

 \begin{conjecture} 
$I(\t (Seg(\BP A_1^*\times\cdots\times \BP A_n^*))$
is generated
the quadrics  in $S^2(A_1\ot \cdots \ot A_m)$ which  have at least
four $\La 2 $ factors, the cubics with  four 
$S_{2,1} $   factors and all other factors $S_{3,0}$, and
the quartics with three $S_{2,2} $'s and all other factors
$S_{4,0}$.
\end{conjecture}

\vspace{1cm}

{\small
\noindent {\sc Joseph M. Landsberg}, 
Department of Mathematics,
  Texas A\&M University,
  Mailstop 3368,
  College Station, TX 77843-3368, USA

\smallskip
 
\noindent {\sc Jerzy Weyman}, 
Department of Mathematics,
  Northeastern University,
  360 Huntington Avenue,
Boston, MA 02115, USA 
}

\smallskip
 
\end{document}